\title{Counting RSA-integers}
\author{Andreas  Decker and Pieter Moree}
\documentclass[12pt]{article}
\usepackage{amssymb, latexsym, amsfonts}
\def\@ptsize{2}
\newtheorem{Thm}{Theorem}

\newtheorem{Lem}{Lemma}

\newtheorem{cor}{Corollary}

\newcommand{\qed}{\hfill $\Box$}

\begin{document}
\date{}
\maketitle
{\def\thefootnote{}
}

\begin{abstract}
\noindent In the RSA cryptosystem integers  of the form $n=p\cdot q$ with $p$ and
$q$ primes of comparable size (`RSA-integers') play an important role. It
is a folklore result of cryptographers that $C_r(x)$, the number of integers $n\le x$ that
are of the form $n=pq$ with $p$ and $q$ primes such
that $p<q<rp$, is for fixed $r>1$ asymptotically equal to $c_rx\log^{-2}x$ for
some constant $c_r>0$. Here we prove this and show that $c_r=2\log r$.
\end{abstract}
\section{Introduction}
Let $C_r (x)$ denote the number of integers $n\le x$ that
are of the form $n=pq$ with $p$ and $q$ primes such
that $p<q<rp$, where $r>1$ is an arbitrary real number.
In this note we establish the following result.
\begin{Thm} As $x$ tends to infiinity, we have
$$C_r (x) = \frac{2 x\log r}{\log^2 x} + O\big(\frac{r\log(er)x}{\log^3 x}\big).$$
\end{Thm}
\begin{cor}
If $r=o(\log x)$, then $C_r(x)\sim {2x\log r\over \log^2 x}$, as $x$ tends to infinity.
\end{cor}
Since $C_r(x)\le x$, the result is only non-trivial if $r=o(\log^3 x/\log \log x)$.\\
\indent It is informative to compare Theorem 1 with a classical one due 
to Landau \cite[205--213]{L}, who in 1909 proved
that $\pi_2(x)$, the number of integers $n\le x$ of the form $n=pq$ with $p$ and $q$
distinct primes satisfies, as $x$ tends to infinity,
$$\pi_2(x)\sim {2x\log \log x\over \log x}.$$
Since then various authors considered the related problem where $n$ consists of 
precisely $k$ primes factors and $k$ is allowed to vary to some extent with $x$.
For a nice survey, see Hildebrand \cite{H}. Note that $C_x(x)=\pi_2(x)$.

\section{Proof of Theorem 1}
Let $\pi(x)$ denote the number of primes not exceeding $x$. All we need regarding
$\pi(x)$ is the
estimate 
\begin{equation}
\label{pietje}
\pi(x)=\int_2^x{dt\over \log t}+
O\Big({x\over \log ^3x}\Big)={x\over \log x}+
{x\over \log^2 x}+O\Big({x\over \log ^3x}\Big).
\end{equation}
The integral in this estimate is usually denoted by Li$(x)$, the logarithmic integral.
Using this estimate one easily infers the following stronger form of the so called
second theorem of Mertens, which one often encounters in the literature with
error term $O(\log^{-1}z)$, see e.g. Tenenbaum \cite[p. 16]{T}. For 
a version of this result with still better error term, see e.g. Landau \cite[201]{L}.
\begin{Lem} 
\label{mertens}
We have
$$\sum_{p\le z}{1\over p}=\log \log z + c_1+O({1\over \log^2 z}),$$
where $c_1$ is a constant.
\end{Lem}
{\it Proof}. Write $\pi(z)={\rm Li}(z)+E(z)$. By (\ref{pietje}) we have
have $E(z)=O(z\log^{-3}z)$.
By Stieltjes integration we find
$$\sum_{p\le z}{1\over p}=\int_{2}^z{d\pi(t)\over t}=\int_2^z{dt\over t\log t}+
\int_2^{z}{dE(t)\over t}=\log\log z+c_2+\int_2^{z}{dE(t)\over t}.$$
On noting that
$$\int_2^z{dE(t)\over t}=c_3+{E(z)\over z}-\int_z^{\infty}{E(t)dt\over t^2},$$
and that 
$$O(\int_z^{\infty}{E(t)dt\over t^2})=O(\int_z^{\infty}{dt\over t\log^3 t})=
O({1\over \log^2 z}),$$
the result follows. \qed\\

\noindent For any prime $p$ we define $f_p(x)$
to be the number of primes $q$ such that $pq\le x$
and $p<q\le rp$. We clearly have
\begin{equation}
\label{e1}
C_r (x)=\sum_{p\le x}f_p(x).
\end{equation}
\begin{Lem}
\label{1}
We have
$$f_p(x)=\cases{\pi(rp)-\pi(p) & if $p\le \sqrt{{x\over r}}$;\cr
\pi({x\over p})-\pi(p) & if $\sqrt{{x\over r}}<p\le \sqrt{x}$;\cr
0 & if $p> \sqrt{x}$.}$$
\end{Lem}
{\it Proof}. Since $p<q$ and $pq\le x$ we infer that $f_p(x)=0$
for $p> \sqrt{x}$. So let us assume that $p\le \sqrt{x}$. Note
that we have $pq\le x$ and $p<q\le rp$ iff
$$p<q\le \min\Big\{rp,{x\over p}\Big\}.$$
We infer that $$f_p(x)=\pi\Big(\min\Big\{rp,{x\over p}\Big\}\Big)-\pi(p).$$
On noting that
$$\min\Big\{rp,{x\over p}\Big\}=\cases{rp & if $p\le \sqrt{x \over r}$;\cr
{x\over p} & if $p>\sqrt{x\over r}$,}$$
the proof is completed. \qed\\

\noindent On combining (\ref{e1}) and Lemma \ref{1} we find that
\begin{equation}
\label{centraal}
C_r (x)=-\sum_{p\le \sqrt{x}}\pi(p)+\sum_{p\le \sqrt{x\over r}}\pi(rp)
+\sum_{\sqrt{x\over r}<p\le \sqrt{x}}\pi\Big({x\over p}\Big).
\end{equation}
The first two sums in the latter expression can be estimated using
Lemma \ref{2}, the third sum using Lemma \ref{3}.
\begin{Lem}
\label{2}
Let $r\ge 1$. We have
\begin{displaymath}
\sum_{p\le z} \pi(rp) = \frac{rz^2}{2 \log^2 z}+O\Big(\frac{r\log(er)z^2}{\log^3 z}\Big).
\end{displaymath}
\end{Lem}
{\it Proof}. First let us assume that $r\le z^{1/4}$.
Comparison of $\pi(rp)$ and $\pi(p)$ yields
\begin{equation}
\label{compare}
\pi(rp)={rp\over \log p}+O\Big({r\log(er)p\over \log^2p}\Big)=r\pi(p)+
O\Big({r\log(er)p\over \log^2p}\Big).
\end{equation}
Since $\sum_{p\le z}p\le \pi(z)z$ we obtain
\begin{equation}
\label{estimate}
O\Big(\sum_{p\le z}{p\over \log^2p}\Big)=O(\sum_{p\le z^{1/3}}p)
+O\Big({1\over \log^{2}z}\sum_{z^{1/3}<p\le z}p\Big)=O\Big({z^2\over \log^3 z}\Big).
\end{equation}
{}From (\ref{compare}) and (\ref{estimate}), we infer that
\begin{equation}
\label{pfoe}
\sum_{p\le z} \pi(rp) =r\sum_{p\le z}\pi(p)+O\Big(\frac{r\log(er)z^2}{\log^3 z}\Big).
\end{equation}
We see that
$$\sum_{p\le z}\pi(p)=\sum_{i\le \pi(z)}i
={1\over 2}\pi(z)(\pi(z)+1),$$
on noting that as $p$ runs over all primes $\le z$,
$\pi(p)$ runs over $1,2,\ldots,\pi(z)$. Hence
\begin{displaymath}
\sum_{p\le z}\pi(p) = {1\over 2} \Big(\frac{z}{\log z}+O\Big(\frac{z}{\log^2 z}
\Big)\Big) \Big(\frac{z}{\log z}+O\Big(\frac{z}{\log^2 z}\Big)+1\Big)=\frac{z^2}{2\log^2 z}+O\Big(\frac{z^2}{\log^3 z}\Big),
\end{displaymath}
which in combination with (\ref{pfoe}) yields the result in case $r\le z^{1/4}$.\\
\indent In case $r>z^{1/4}$ it is enough to show that
$\sum_{p\le z}\pi(rp)=O(rz^2\log^{-2}z)$. On noting that $\pi(rp)=O(rp/\log p)$, this
estimate easily follows (analogous to the derivation of (\ref{estimate})).\qed

\begin{Lem}
\label{3}
Suppose that $1\le r\le \sqrt{x}$. Then
$$\sum_{\sqrt{x\over r}<p\le \sqrt{x}}\pi\Big({x\over p}\Big)={2x\log r\over \log^2 x}+
O\Big({x\log^2 (er)\over \log^3 x}\Big).$$
\end{Lem}
{\it Proof}. On writing $p=\sqrt{x\over a}$, we find that, for $\sqrt{x\over r}<p\le x$, we
have
\begin{equation}
\label{H}
\pi({x\over p})={2x\over p\log x}+O\Big({x\log (er)\over p\log^2x}\Big).
\end{equation}
We infer from Lemma \ref{mertens} that
$$H(x):=\sum_{\sqrt{x\over r}<p\le \sqrt{x}}{1\over p}=-\log\Big(1-{\log r\over \log x}\Big)
+O({1\over \log^2(x/r)}\Big)={\log r\over \log x}+O\Big({\log^2(er)\over \log^2x}\Big).$$
{}From (\ref{H}) we infer that 
$$\sum_{\sqrt{x\over r}<p\le \sqrt{x}}\pi\Big({x\over p}\Big)=2H(x){x\over \log x}
+O\Big(H(x){x\log (er)\over \log^2x}\Big).$$
This, together with the estimate for $H(x)$ we determined, yields the result. \qed\\

\noindent {\it Proof of Theorem 1}. 
For $r>\sqrt{x}$ we have $C_r(x)\le x=O(r\log(er)x\log^{-3}x)$ and so
Theorem 1 is trivially true for this $r$-range. Now assume that $r<\sqrt{x}$.
By Lemma \ref{2} (with $r=1$ and $z=\sqrt{x}$) we have
$$\sum_{p\le \sqrt{x}}\pi(p)={2x\over \log^2x}+O\Big({x\over \log^3 x}\Big).$$
Lemma \ref{2} with $z=\sqrt{{x\over r}}$ yields 
$$\sum_{p\le \sqrt{{x\over r}}}\pi(rp)={2x\over \log^2x}+O\Big({xr\log(er)\over \log^3 x}\Big).$$
The proof now follows on inserting the latter two estimates and the estimate
given in Lemma \ref{3} in the equality (\ref{centraal}). \qed\\

\noindent {\bf Acknowledgement}. (A.D.:) This note was written whilst the 
first author did an internship at the Max-Planck-Institut f\"ur Mathematik in Bonn. 
The other interns Alexander Bridi, Patrizia Dressler, Silke Glas and Thorge Jensen also 
contributed to this paper with different suggestions for improvement and some own 
calculations. Last, but not least, the results on this problem are also due to the 
good understanding of the interns, their good relationship with their 
coordinator, the second author, and the pleasant atmosphere at the MPIM in Bonn.\\
\indent  (P.M.:) The problem of estimating
RSA-integers was proposed by Benne de Weger (TU Eindhoven) to the second author, who
gave it as one of several problems to the above interns. It was in essence solved
by Andreas Decker. The proof presented here is shorter and somewhat different and
has absolute implicit errors instead of $r$-dependent ones.

\vfil\eject

\medskip\noindent {\footnotesize Achtern Diek 32, D-49377 Vechta, Germany.\\
e-mail: {\tt andreasd@uni-bonn.de}}\\

\medskip\noindent {\footnotesize Max-Planck-Institut f\"ur Mathematik,\\
Vivatsgasse 7, D-53111 Bonn, Germany.\\
e-mail: {\tt moree@mpim-bonn.mpg.de}\\
(corresponding author)}
\vskip 5mm
\end{document}